\documentclass[draft,preprint,12pt,numbers,sort&compress]{elsarticle}
\usepackage{}
\usepackage{mathrsfs}
\usepackage{amssymb}
\usepackage{amsthm}
\usepackage{mathrsfs}
\usepackage[centertags]{amsmath}
\usepackage{amsfonts}

\usepackage{color}

\input amssym.def
\input amssym.tex

\date{}

\newtheorem{Theorem}{Theorem}[section]

\newtheorem{Lemma}{Lemma}[section]

\newcommand\R{\mbox{\bf R}}
\newcommand\T{\mbox{\bf T}}

\newcommand\Z{\mbox{\bf Z}}
\newcommand\z{\mbox{\bf z}}
\newcommand\SR{\mbox{\scriptsize\bf R}}

\newcommand{\definition}{{\lower .5ex
  \hbox{$\>\>\stackrel{\triangle}{=}\>\>$} }}
\newcommand\supp{\mathop{\rm supp}}

\begin{document}

\baselineskip=22pt
\thispagestyle{empty}

\begin{center}
{\Large \bf Convergence problem of Ostrovsky }\\[1ex]
{\Large \bf equation with rough data and random data}\\[1ex]

{\quad Wei Yan\footnote{Email:  011133@htu.edu.cn}$^{a}$, Qiaoqiao Zhang$^{a}$,
Jinqiao Duan\footnote{Email:  duan@iit.edu}$^{b}$, Meihua Yang\footnote{ Email: yangmeih@hust.edu.cn}$^{c*}$}\\[1ex]

{$^a$School of Mathematics and Information Science, Henan
Normal University,}\\
{Xinxiang, Henan 453007,   China}\\[1ex]

{$^b$Department of Applied Mathematics, Illinois Institute of Technology,}\\[1ex]
{Chicago, IL 60616, USA}\\[1ex]

{$^c$School of Mathematics and Statistics, Huazhong University of Science and Technology, }\\
{Wuhan, Hubei 430074,  China}\\[1ex]

\end{center}
\noindent{\bf Abstract.} In this paper, we consider  the convergence
 problem of free Ostrovsky  equation with rough data and random data respectively.
 We show
 the almost everywhere pointwise convergence
of free Ostrovsky equation
with initial rough data  in $H^{s}(\R)$ for $s\geq \frac{1}{4}$. Counterexample
 is constructed to show that the maximal function estimate related to the free Ostrovsky equation  can
 fail if $s<\frac{1}{4}$.
We also show
 the stochastic continuity at $t=0$
of free Ostrovsky equation  with initial
random data in $L^{2}(\R)$.

 \bigskip

\noindent {\bf Keywords}: Stochastic  pointwise convergence;
Free Ostrovsky  equation; Rough data; Random data
\medskip

\medskip
\noindent {\bf Corresponding Author:}Meihua Yang

\medskip
\noindent {\bf Email Address:}yangmeih@hust.edu.cn

\medskip
\noindent {\bf AMS  Subject Classification}:   42B25; 42B15; 35Q53

\leftskip 0 true cm \rightskip 0 true cm

\newpage

\baselineskip=20pt

\bigskip
\bigskip

{\large\bf 1. Introduction}
\bigskip

\setcounter{Theorem}{0} \setcounter{Lemma}{0}\setcounter{Definition}{0}\setcounter{Proposition}{1}

\setcounter{section}{1}

In this paper,  we  investigate  the pointwise convergence  problem of  the free Ostrovsky equation
\begin{eqnarray}
&&u_t+\partial_{x}^{3}u\pm \partial_{x}^{-1}u=0, \label{1.01}\\
&&u(x,0)=f(x)   \label{1.02}.
\end{eqnarray}
It is easily checked that $U(t)f=\frac{1}{\sqrt{2\pi}}\int_{\SR} e^{ix\xi+it(\xi^{3}\pm\frac{1}{\xi})}
\mathscr{F}_{x}{f}(\xi)d\xi$ is the solution to (\ref{1.01})-(\ref{1.02}),  where $\mathscr{F}_{x}$ is the
Fourier transform in the $x$ variable, see the concrete definition of $\mathscr{F}_{x}$ in the line 13 of Page 5.
Ostrovsky equation
was proposed by Ostrovsky  \cite{GiGrSt,Gri,O} as a model
for weakly nonlinear long waves in a rotating liquid, by taking
into account of the Coriolis force.
It describes the propagation of surface waves
in the ocean in a
rotating frame of reference.  Cauchy problems for the Ostrovsky equation are
 investigated in
 \cite{CR,GHO,GL,LHY,HJ0,
IM1,IM2,IM3,IM4,IM5,LM2006,Tsu,VL,YLHD}.

Carleson \cite{Carleson} initiated the
pointwise converge problem, more precisely, Carleson showed pointwise
convergence problem of the one dimensional Schr\"odinger equation in
$H^{s}(\R)$ with $s\geq \frac{1}{4}$. Dahlberg and Kenig \cite{DK} showed
that the pointwise convergence of the
Schr\"odinger equation does not hold for
$s<\frac{1}{4}$ in any dimension.  Dahlberg and Kenig \cite{DK} and Kenig et al. \cite{KPV1991,KPV1993}
have proved that the pointwise convergence of  KdV  equation  holds if and only if $s\geq \frac{1}{4}$.

For the  pointwise  convergence  problem of
 Schr\"odinger  equations  in  higher dimension,  Bourgain\cite{Bourgain2016}
 recently presented counterexamples showing that convergence of Schr\"odinger equation
 in $\mathbb{R}^n$ can fail if $s<\frac{n}{2(n+1)},n\geq2$.
 Du et al. \cite{DGL} proved that the pointwise convergence problem of two dimensional
  Schr\"odinger
 equation in $H^{s}(\R^{2})$ with $s>\frac{1}{3}$.
 Du and Zhang \cite{DZ} proved that the
  pointwise convergence problem of $n$ dimensional Schr\"odinger
 equation in $H^{s}(\R^{n})$ with $s>\frac{n}{2(n+1)}, n\geq3.$ See more references such as
  \cite{Bourgain1995,Bourgain2013,CLV,Cowling,DG, DGZ,GS,Lee,LR2015,LR2017,MVV,S,Shao,Tao,TV,Vega,Zhang}.
  Miao et al. \cite{MYZ2015,MZZ2015} studied the  pointwise  convergence  problem of 2D
 fractional order Schr\"odinger equations  and  Schr\"odinger
 equation with inverse-square potential, respectively. The pointwise
 convergence problem of Schr\"odinger equation on the torus $\mathrm{\T}^{n}$ was tackled first
 by Moyua-Vega \cite{MV} and recently extended by Wang-Zhang \cite{WZ},  Compaan-Luc$\grave{a}$-Staffilani \cite{CLS} and
 Eceizabarrena-Luc$\grave{a}$\cite{EL}, etc.

The method of applying suitable randomized initial data was first introduced by
Lebowitz-Rose-Speer  in \cite{LRS} and Bourgain \cite{Bourgain1994, Bourgain1996}
 and Burq-Tzvetkov \cite{BTL,BTG}.
This method was also applied to study nonlinear dispersive equations and
 hyperbolic equations in scaling super-critical regimes, for example,
see \cite{BOP,BOP2015,BOP2017,CO,CG,CZ,DC,DG,HM,HK,KM,LMCPDE,M,NORS,NPS,NS,OOP,OP,P,ZF2011,ZF2012}.
Very recently, Compaan et al.
 \cite{CLS} applied randomized  initial data to study
  pointwise convergence  of the Schr\"odinger flow.

In this paper, motivated by \cite{CLS,KPV1991,KPV1993,Du},
we investigate the convergence problem
of free Ostrovsky equation with rough data and random data.
We show
 the almost everywhere pointwise convergence
of free Ostrovsky equation
with initial rough data  in $H^{s}(\R)$ for $s\geq \frac{1}{4}$. Counterexample is
 constructed to show that the maximal function estimate related to the
 free Ostrovsky equation  can
 fail if $s<\frac{1}{4}$.
We also show
 the stochastic continuity
of free Ostrovsky equation  with initial
random data in $L^{2}(\R)$.
The main ingredients are  the
 density theorem, high-low frequency idea,  Wiener  decomposition of frequency spaces
 and Lemmas 2.1-2.7 as well as some probabilistic estimates.
The main difficulty is that zero is  the singular point
of the phase functions $\xi^{3}\pm \frac{1}{\xi}$ of free Ostrovsky
equation.

Now we present the deterministic results concerning the pointwise convergence problem which are just Theorems 1.1, 1.2.

\begin{Theorem} \label{Theorem1}(Pointwise convergence)
Let  $f\in  H^{s}(\R)$ with $s\geq \frac{1}{4}.$ Then, we have
\begin{eqnarray}
\lim\limits_{t\longrightarrow 0}U(t)f(x)=f(x)\label{1.09}
\end{eqnarray}
almost everywhere with respect to $x$.

\end{Theorem}

\begin{Theorem}[] \label{Theorem2}
For $s<\frac{1}{4}$ and $f_{k}=\frac{1}{\sqrt{2\pi}}\int_{\SR}e^{ix\xi}2^{-k(s+\frac{1}{2})}
\chi_{2^{k}\leq |\xi|\leq 2^{k+1}}(\xi)d\xi$, we have
\begin{align*}
\lim\limits_{k\rightarrow \infty}\frac{\left\|\sup\limits_{t>0}|U(t)f_{k}|\right\|_{L_{x}^{4}}}{\|f_{k}\|_{H^{s}(\SR)}}=\infty.
\end{align*}
The  maximal  inequality
\begin{eqnarray}
\|U(t)f\|_{L_{x}^{4}L_{t}^{\infty}}
\leq C\|f\|_{H^{s}(\SR)}\label{1.010}
\end{eqnarray}
does not hold in general for $f\in H^{s}(\R)$
if $s<\frac{1}{4}.$
\end{Theorem}

Now we introduce the randomization procedure for the initial data, which can be
 seen in \cite{BOP,BOP2015,LMCPDE,ZF2012}.
 Let $\psi\in\mathcal{S}(\R)$ be an even, non-negative jump function
with $supp(\psi)\subseteq [-1,1]$ and such that for $\xi\in\R$,
 \begin{eqnarray}
\sum_{k\in
\z}\psi(\xi-k)=1.\label{1.04}
\end{eqnarray}
For every $k\in
\Z$, we define the function $\psi(D-k)f:\R\rightarrow\mathbb{C}$ by
\begin{eqnarray*}
(\psi(D-k)f)(x)=\mathscr{F}^{-1}_{\xi}\big(\psi(\xi-k)\mathscr{F}_xf\big)(x), ~x\in \R.
\end{eqnarray*}
Note that these projections satisfy
 a unit-scale Bernstein inequality  which can be seen in Lemma 2.1 of \cite{LMCPDE},
 namely that,
 for any $p_1,p_2$, which satisfies that
$2 \leq p_1 \leq p_2 \leq\infty$,
there exists a $C= C(p_1, p_2)>0$ such that for any $f\in L^2(\R)$
and $k\in \Z$,
 \begin{eqnarray}
\left\|\psi(D-k)f\right\|_{L_{x}^{p_2}(\SR)}\leq C
\left\|\psi(D-k)f\right\|_{L_{x}^{p_1}(\SR)}\leq C
\left\|\psi(D-k)f\right\|_{L_{x}^{2}(\SR)}.\label{1.05}
\end{eqnarray}

Let $\{g_k\}_{k\in Z}$ be a sequence of independent, zero-mean, complex-valued
 Gaussian random variables
 on a probability space $(\Omega,\mathcal{A}, \mathbb{P})$,
 where the real and imaginary parts of $g_k$ are independent
  and endowed with probability
distributions $\mu_k^1$ and $\mu_k^2$ respectively. The probability
distributions $\mu_k^1$ and $\mu_k^2$ satisfy the following condition:

There exists $c>0$ such that
 \begin{eqnarray}
\Big|\int_{-\infty}^{+\infty}e^{\gamma x}d\mu_k^j(x)\Big|\leq e^{c\gamma^2},\quad\text{for all}~ \gamma\in \R, k\in  \Z, j=1, 2.\label{1.06}
\end{eqnarray}
Note that, since the real and imaginary parts of $g_k$ possesses the same density function $\frac{1}{\sqrt{2\pi}}e^{-\frac{x^{2}}{2}},$ thus, we have
  \begin{eqnarray*}
&&\int_{-\infty}^{+\infty}e^{\gamma x}d\mu_k^j(x)=\frac{1}{\sqrt{2\pi}}\int_{\SR}e^{-\frac{x^{2}}{2}+\gamma x}dx\nonumber\\
&&=e^{\frac{\gamma^{2}}{2}}\frac{1}{\sqrt{2\pi}}\int_{\SR}e^{-\frac{(x-\gamma)^{2}}{2}}dx\nonumber\\&&=e^{\frac{\gamma^{2}}{2}}
,\quad\text{for all}~ \gamma\in \R, k\in  \Z, j=1, 2.
\end{eqnarray*}
Which means property \eqref{1.06} is satisfied.

Thereafter for a given $f\in H^{s}(\R)$, we define its
randomization by
\begin{eqnarray}
f^\omega:=\sum_{k\in \z}g_k(\omega)\psi(D-k)f.\label{1.07}
\end{eqnarray}
We define
\begin{eqnarray*}
\|f\|_{L_{\omega}^{p}(\Omega)}=\left[\int_{\Omega}|f(\omega)|^{p}d\mathbb{P}(\omega)\right]^{\frac{1}{p}}.
\end{eqnarray*}

If $f\in H^{s}(\R)$,
then the randomized function $f^\omega$
is almost surely
in $H^{s}(\R)$ and $\|\|f^\omega\|_{H^{s}}\|_{L_{\omega}^{2}}=\|f\|_{H^{s}}$,
see Lemma 2.2 in \cite{BOP}. This randomization improves the
 integrability of $f$, see Lemma 2.3 of \cite{BOP}. Such results for random
 Fourier series are known as Paley-Zygmund's theorem \cite{PZ}.

 We will restrict ourselves to a subset $\Sigma\subset \Omega$ with
$P(\Sigma)=1$ such that $f^{\omega}\in H^{s}$ for every $\omega \in \Sigma.$

Now we state Theorem 1.3 as follows.

\begin{Theorem} \label{Theorem3}(Stochastic continuity )
Let $f\in  L^{2}(\R)$  and $f^{\omega}$ be a randomization of $f$ as
defined in (\ref{1.07}). Then,  $\forall \alpha>0$, we have
\begin{eqnarray}
\lim\limits_{t\longrightarrow 0}\mathbb{P}\left(\omega\in\Omega:|U(t)f^{\omega}(x)- f^{\omega}(x)|>\alpha\right)=0,
\end{eqnarray}
which is independent of   $x$. More precisely, $\forall \epsilon>0$ such that $2Ce\epsilon(\ln\frac{3C_{1}}{\epsilon})^{\frac{1}{2}}<\alpha$ and   when  $|t|<\epsilon^{2}$, we have
\begin{eqnarray}
&&\mathbb{P}\left(\left\{\omega\in \Omega:
\left|U(t)f^{\omega}-f^{\omega}\right|>\alpha\right\}\right)\leq\epsilon.\label{1.011}
\end{eqnarray}
Here, $C,C_{1}$ appear in Lemmas 3.2-3.4.
\end{Theorem}

\noindent{\bf Remark 1.} The definition of stochastic continuity was given in page 70 of \cite{DPZ}.

Now we present some notations.
\begin{eqnarray*}
&&\mathscr{F}_{x}f(\xi)=\frac{1}{\sqrt{2\pi}}\int_{\SR}e^{-ix \xi}f(x)dx,\\
&&\mathscr{F}_{x}^{-1}f(\xi)=\frac{1}{\sqrt{2\pi}}\int_{\SR}e^{ix \xi}f(x)dx,\\
&&\|f\|_{L_{x}^{q}L_{t}^{p}}=\left(\int_{\SR}
\left(\int_{\SR}|f(x,t)|^{p}dt\right)^{\frac{q}{p}}dx\right)^{\frac{1}{q}},\\
&&\|f\|_{L_{xt}^{p}}=\|f\|_{L_{x}^{p}L_{t}^{p}}.
\end{eqnarray*}

$H^{s}(\R)=\left\{f\in \mathscr{S}^{'}(\R):\|f \|_{H^{s}(\SR)}=
 \|\langle\xi\rangle ^{s}\mathscr{F}_{x}{f}\|_{L_{\xi}^{2}(\SR)}<\infty\right\},$
 where $\langle\xi\rangle ^{s} =(1+\xi^{2})^{\frac{s}{2}}$ for any $\xi\in
\R$.

$|E|$ denotes by the Lebesgue measure of set $E$.

Let $\phi$ be a smooth jump function such that $\phi(\xi)=1$
 for $|\xi|\leq 1$ and $\phi(\xi)=0$  for $|\xi|>2$. Then, we define
for every dyadic integer $N\in 2^{\z}$,
\begin{eqnarray*}
&&\mathscr{F}_{x}P_{N}f(\xi)=\left[\phi\left(\frac{\xi}{N}\right)-
\phi\left(\frac{2\xi}{N}\right)\right]\mathscr{F}_{x}f(\xi),\\
&&\mathscr{F}_{x}P_{\leq N}f(\xi)=\phi\left(\frac{\xi}{N}\right)
\mathscr{F}_{x}f(\xi),\\
&&\mathscr{F}_{x}P_{>N}f(\xi)=\left[1-
\phi\left(\frac{\xi}{N}\right)\right]\mathscr{F}_{x}f(\xi).
\end{eqnarray*}

\noindent{\bf Remark 2.} Now we give the outline of proof of Theorem 1.1.
For the Ostrovsky equation which possesses the phase function
 $\xi^{3}\pm \frac{1}{\xi},$ the main difficulty is that zero is  the singular point
of the phase functions.

In order to deal with the singular point, we use the high-low frequency idea,
that is, we establish estimates for high frequency,
$|\xi^{3}\pm \frac{1}{\xi}|\sim |\xi|^{3}$ for $|\xi|\geq 8$,
and low frequency, $|\xi^{3}\pm \frac{1}{\xi}|\simeq |\xi|^{-1}$ for $|\xi|\leq 8$, separately.

Note that,  by density theorem which is just  Lemma 2.2 in \cite{Du},
for any $f\in H^{s}(\R)$ with $s\geq \frac{1}{4},$  the following decomposing properties
hold:

$\forall\epsilon>0$, $f$ can be
decomposed as $f=g+h$, where $g$ is a rapidly decreasing function,
 $\|h\|_{H^{s}(\SR)}<\epsilon.$

Hence, we can establish the corresponding estimates.

Concretely, on one hand,
for the high frequency: $|\xi^{3}\pm \frac{1}{\xi}|\sim |\xi|^{3}$ for $|\xi|\geq 8$, since $g$
 is a rapidly decreasing function, following the method of Lemma 2.3 in \cite{Du}, we prove
\begin{eqnarray}
\left|U(t)P_{\geq8}g-P_{\geq8}g\right|\longrightarrow 0,\quad\text{as}~ t\rightarrow 0. \label{1.012}
\end{eqnarray}
For the detail of \eqref{1.012}, we refer the
readers to Lemma 2.3 in this paper.

And for the low frequency,  since $g$ is a rapidly decreasing function,  with the aid of Lemma 2.4 obtained in this paper,
we establish
\begin{eqnarray}
\left|U(t)P_{\leq8}g-P_{\leq8}g\right|\longrightarrow 0,\quad\text{as}~ t\rightarrow 0. \label{1.013}
\end{eqnarray}
On the other hand, by using Lemma 2.1, we have
\begin{eqnarray}
\left\|U(t)P_{\geq8}h\right\|_{L_{x}^{4}L_{t}^{\infty}}\leq C\|h\|_{H^{\frac{1}{4}}(\SR)}\label{1.014}
\end{eqnarray}
with $\|h\|_{H^{\frac{1}{4}}(\SR)}<\epsilon$.
Since   $\epsilon$  can be  chosen  as small as needed, following the method of Lemma 2.3 in \cite{Du},  from (\ref{1.014}),   we obtained
\begin{eqnarray}
\left|U(t)P_{\geq8}h-P_{\geq8}h\right|\longrightarrow 0,\quad\text{as}~ t\rightarrow 0.\label{1.016}
\end{eqnarray}
Furthermore, using Lemma 2.2 in this  paper,  we obtain that there exists $\delta_\epsilon>0$
such that when $|t|\leq\frac{\delta_\epsilon}{C}$,
\begin{eqnarray}
\left|U(t)P_{\leq8}h-P_{\leq8}h\right|\leq 2\epsilon\label{1.018}.
\end{eqnarray}
with $\|h\|_{H^{\frac{1}{4}}}<\epsilon$.
From (\ref{1.012}), \eqref{1.013},  \eqref{1.016} and (\ref{1.018}), as $t\longrightarrow 0,$
we have
\begin{eqnarray}
U(t)f\longrightarrow f\label{1.019}.
\end{eqnarray}

\noindent{\bf Remark 3.}
By presenting particular initial data, we give a counterexample to show that the maximal function
 estimate can be invalid for $s<\frac{1}{4}.$ Then, we obtained Theorem 1.2.

\noindent{\bf Remark 4.}
Now, we present the proof of Theorem 1.3.

By density theorem, that is, rapidly decreasing functions are dense in $L^2(\R)$,
for any $f\in L^{2}(\R)$, the following decomposing properties
hold:

$\forall\epsilon>0$, $f$ can be
decomposed as $f=g+h$,  where $g$ is a rapidly
 decreasing function and  $\|h\|_{L^{2}(\SR)}<\epsilon$.

Then,
\begin{align*}
f^{\omega}=g^{\omega}+h^{\omega}
\end{align*}
and
\begin{eqnarray}
&&U(t)f^{\omega}-f^{\omega}=U(t)g^{\omega}-g^{\omega}
+U(t)h^{\omega}-h^{\omega}.\label{1.020}
\end{eqnarray}
Here, $f^{\omega},g^{\omega}$ and $h^{\omega}$ are defined as in (\ref{1.07}).

 $\forall\alpha>0$, by using Lemma 3.1 and high-low frequency technique, since $g$ is a rapidly
 decreasing function, $\forall\epsilon>0$ and $|t|>0$, we have
 \begin{eqnarray}
&&\mathbb{P}\left(\left\{\omega\in \Omega:
\left|U(t)g^{\omega}-g^{\omega}\right|>\frac{\alpha}{2}\right\}\right)\leq C_{1} e^{-\left(\frac{\alpha}
{Ce\left[\epsilon+\frac{|t|}{\epsilon}\right]}\right)^{2}}.\label{1.021}
\end{eqnarray}
For the details of proof, we refer the readers to Lemma 3.2.

By using
Lemmas 3.1,  2.7, we have
\begin{eqnarray}
\mathbb{P}\left(\left\{\omega\in \Omega:
|U(t)h^{\omega}|>\frac{\alpha}{4}\right\}\right)\leq C_{1}e^{-\left(\frac{\alpha}{Ce\|h\|_{L^{2}}}\right)^{2}}.\label{1.022}
\end{eqnarray}
For the details of proof, we refer the readers to Lemma 3.3.

By using  Lemmas 3.1, 2.6, we have
\begin{eqnarray}
\mathbb{P}\left(\left\{\omega\in \Omega: |h^{\omega}|>\frac{\alpha}{4}\right\}\right)\leq C_{1}e^{-\left(\frac{\alpha}{Ce\|h\|_{L^{2}}}\right)^{2}}\label{1.023}.
\end{eqnarray}
Thus, combining (\ref{1.020}) with (\ref{1.021})-(\ref{1.023}), $\forall \epsilon>0$
such that $2Ce\epsilon(\ln\frac{3C_{1}}{\epsilon})^{\frac{1}{2}}<\alpha$ and taking $\|h\|_{L^{2}}\leq \epsilon$, when  $|t|<\epsilon^{2}$
\begin{eqnarray}
&&\mathbb{P}\left(\left\{\omega\in \Omega:
\left|U(t)f^{\omega}-f^{\omega}\right|>\alpha\right\}\right)\nonumber\\&&\leq
\mathbb{P}\left(\left\{\omega\in \Omega:
\left|U(t)g^{\omega}-g^{\omega}\right|>\frac{\alpha}{2}\right\}\right)+\mathbb{P}\left(\left\{\omega\in \Omega:
|U(t)h^{\omega}|>\frac{\alpha}{4}\right\}\right)\nonumber\\&&\qquad+\mathbb{P}\left(\left\{\omega\in \Omega: |h^{\omega}|>\frac{\alpha}{4}\right\}\right)
\nonumber\\&&\leq  C_{1} e^{-\left(\frac{\alpha}
{Ce\left[\epsilon+\frac{|t|}{\epsilon}\right]}\right)^{2}}+2C_{1}e^{-\left(\frac{\alpha}{Ce\|h\|_{L^{2}}}\right)^{2}}\nonumber\\&&\leq  C_{1} e^{-\left(\frac{\alpha}
{Ce\left[\epsilon+\frac{|t|}{\epsilon}\right]}\right)^{2}}+2C_{1}e^{-\left(\frac{\alpha}{Ce\epsilon}\right)^{2}}\nonumber\\
&&\leq C_{1} e^{-\left(\frac{\alpha}{2Ce\epsilon}\right)^{2}}+2C_{1}e^{-\left(\frac{\alpha}{Ce\epsilon}\right)^{2}} \nonumber\\
&&\leq 3C_{1}e^{-\left(\frac{\alpha}{2Ce\epsilon}\right)^{2}}\leq\epsilon.\label{1.024}
\end{eqnarray}
Hence, for any $\alpha>0$, we have
\begin{eqnarray}
\lim\limits_{t\longrightarrow 0}\mathbb{P}\left(\left\{\omega\in \Omega:
\left|U(t)f^{\omega}-f^{\omega}\right|>\alpha\right\}\right)=0.\label{1.025}
\end{eqnarray}
uniformly with respect to $x$.

This completes the proof of Theorem 1.3.

\bigskip

\bigskip

\setcounter{section}{2}

\noindent{\large\bf 2. Preliminaries }

\setcounter{equation}{0}

\setcounter{Theorem}{0}

\setcounter{Lemma}{0}

\setcounter{section}{2}
In this section, we present some preliminary estimates
related to  Ostrovsky equation. More precisely, Lemmas 2.1-2.4 are used to establish Theorem 1.1 and
Lemmas 2.5-2.7 are used to establish Theorem 1.3.

\begin{Lemma}\label{lem2.1}(Maximal function estimate related to Ostrovsky  equation)
For $f\in H^{\frac{1}{4}}(\R)$, we have
\begin{eqnarray}
&&\left\|U(t)P_{\geq8}f\right\|_{L_{x}^{4}L_{t}^{\infty}}\leq C
\left\|f\right\|_{H^{\frac{1}{4}}(\SR)}\label{2.01}.
\end{eqnarray}
\end{Lemma}

For the proof of Lemma 2.1, we refer the readers to (2.2)  of \cite{GHO}.

\begin{Lemma}\label{lem2.2}(Estimate related to Ostrovsky equation with low frequency)
  $\forall \epsilon>0$ and  $g\in L^{2}(\R)$,
 there exists $\delta_\epsilon>0$ such that
\begin{align}
\left|U(t)P_{\leq8}g-P_{\leq 8}g\right|
\leq \epsilon+\frac{C|t|}{\delta_\epsilon}\|g\|_{L^{2}(\SR)}.\label{2.02}
\end{align}
Here, $C$ is independent of $g,x,t,\epsilon.$

In particular, when $|t|\leq \frac{\delta_\epsilon}{C},\|g\|_{L^{2}}<\epsilon,$ we have
\begin{eqnarray*}
\left|U(t)P_{\leq8}g-P_{\leq 8}g\right|
\leq 2\epsilon.
\end{eqnarray*}
\end{Lemma}
\noindent{\bf Proof.}   $\forall \epsilon>0$, since $g\in L^{2}(\R),$ there exists $\delta_\epsilon >0(<\frac{1}{2})$ such that
\begin{eqnarray}
\left[\int_{|\xi|\leq \delta_{\epsilon}}
|\mathscr{F}_{x}g(\xi)|^{2}d\xi\right]^{\frac{1}{2}}\leq\epsilon.\label{2.03}
\end{eqnarray}
By using the Cauchy-Schwarz inequality and (\ref{2.03}), we have
\begin{eqnarray}
&&\int_{|\xi|\leq \delta_{\epsilon}}|\mathscr{F}_{x}g(\xi)|d\xi\leq
\left[\int_{|\xi|\leq \delta_{\epsilon}}|\mathscr{F}_{x}g(\xi)|^{2}d\xi\right]^{\frac{1}{2}}(2\delta_{\epsilon})^{\frac{1}{2}}\leq
\epsilon.\label{2.04}
\end{eqnarray}
By using the Cauchy-Schwarz inequality, we have
\begin{eqnarray}
&&\int_{\delta_{\epsilon}\leq |\xi|\leq 8}|\mathscr{F}_{x}g(\xi)|d\xi\leq
\left[\int_{\delta_{\epsilon}\leq |\xi|\leq 8}|\mathscr{F}_{x}
g(\xi)|^{2}d\xi\right]^{\frac{1}{2}}\left[\int_{\delta_{\epsilon}\leq |\xi|\leq 8}
d\xi\right]^{\frac{1}{2}}\nonumber\\
&&\leq 3\left[\int_{\delta_{\epsilon}\leq |\xi|\leq 8}
|\mathscr{F}_{x}g(\xi)|^{2}d\xi\right]^{\frac{1}{2}}
\leq 3\|g\|_{L^{2}}.\label{2.05}
\end{eqnarray}
For $\delta_{\epsilon}\leq |\xi|\leq 8$, we have
\begin{eqnarray}
&&\left|e^{it(\xi^{3}\pm \frac{1}{\xi})}-1\right|\leq |t|\left|\xi^{3}
\pm \frac{1}{\xi}\right|\leq \frac{C|t|}{|\xi|}\leq \frac{C|t|}{\delta_{\epsilon}}\label{2.06}.
\end{eqnarray}
Thus, from (\ref{2.03})-(\ref{2.06}), we have
\begin{eqnarray}
&&\left|U(t)P_{\leq 8}g-P_{\leq 8}g\right|
=\left|\int_{|\xi|\leq8}e^{ix \xi}\left[e^{it(\xi^{3}\pm \frac{1}{\xi})}-1\right]
\mathscr{F}_{x}g(\xi)d\xi\right| \nonumber\\
&&\leq \left|\int_{|\xi|\leq \delta_{\epsilon}}e^{ix \xi}
\left[e^{it(\xi^{3}\pm \frac{1}{\xi})}-1\right]
\mathscr{F}_{x}g(\xi)d\xi\right|
\nonumber\\&&\qquad\qquad+\left|\int_{\delta_{\epsilon}
\leq |\xi|\leq 8}e^{ix \xi}\left[e^{it(\xi^{3}\pm
 \frac{1}{\xi})}-1\right]\mathscr{F}_{x}g(\xi)d\xi\right|\nonumber\\
&&\leq \int_{|\xi|\leq \delta_{\epsilon}}|\mathscr{F}_{x}g(\xi)
|d\xi +C|t|\int_{\delta_{\epsilon}\leq |\xi|\leq 8}
\frac{1}{|\xi|}|\mathscr{F}_{x}g(\xi)|d\xi\nonumber\\
&&\leq \epsilon+\frac{C|t|}{\delta_{\epsilon}}\int_{|\xi|\leq8}|\mathscr{F}_{x}g(\xi)
|d\xi\leq \epsilon+\frac{C|t|}{\delta_{\epsilon}}\|g\|_{L^{2}}.\label{2.07}
\end{eqnarray}

This completes the proof of Lemma 2.2.

\begin{Lemma}\label{lem2.3}(Estimate related to Ostrovsky equation with high frequency)
Let $g$ be a rapidly decreasing function. Then, we have
 \begin{align}
\left|U(t)P_{\geq8}g-P_{\geq8}g\right|
\leq C|t|.\label{2.08}
\end{align}
Here, $C$ is independent of $x,t.$
\end{Lemma}
\noindent{\bf Proof.}
Since $g$ is a rapidly decreasing function,  we have
\begin{eqnarray}
\left|U(t)P_{\geq8}g-P_{\geq8}g\right|
&&\leq \left|\int_{|\xi|\geq 8}e^{ix \xi}
\left[e^{it(\xi^{3}\pm \frac{1}{\xi})}-1\right]\mathscr{F}_{x}g(\xi)d\xi\right|
\nonumber\\
&&\leq C|t|\int_{|\xi|\geq8}\left|\xi^{3}\pm \frac{1}{\xi}\right||\mathscr{F}_{x}g(\xi)|d\xi
\nonumber\\
&&\leq C|t|\int_{|\xi|\geq8}\left|\xi\right|^{3}|\mathscr{F}_{x}g(\xi)|d\xi\leq C|t|.\label{2.09}
\end{eqnarray}

This completes the proof of Lemma 2.3.

\begin{Lemma}\label{lem2.4}(Estimate related to Ostrovsky equation with low frequency)
   Let $g$ be a  rapidly decreasing  function, $\forall \epsilon>0(<10^{-2})$,  we have
\begin{align}
\left|U(t)P_{\leq8}g-P_{\leq 8}g\right|
\leq C\left[\epsilon+\frac{|t|}{\epsilon}\right].\label{2.010}
\end{align}
Here, $C$ is independent of $\epsilon,x,t.$
\end{Lemma}
\noindent{\bf Proof.}   $\forall \epsilon>0$, since $g$ is a  rapidly decreasing  function, we have
\begin{eqnarray}
&&\int_{|\xi|\leq \epsilon}|\mathscr{F}_{x}g(\xi)|d\xi\leq C
\epsilon.\label{2.011}
\end{eqnarray}
For $\epsilon\leq |\xi|\leq 8$, we have
\begin{eqnarray}
&&\left|e^{it(\xi^{3}\pm \frac{1}{\xi})}-1\right|\leq |t|\left|\xi^{3}
\pm \frac{1}{\xi}\right|\leq \frac{C|t|}{|\xi|}\leq \frac{C|t|}{\epsilon}\label{2.012}.
\end{eqnarray}
Thus, from (\ref{2.011})-(\ref{2.012}), we have
\begin{eqnarray}
&&\left|U(t)P_{\leq 8}g-P_{\leq 8}g\right|
=\left|\int_{|\xi|\leq8}e^{ix \xi}\left[e^{it(\xi^{3}\pm \frac{1}{\xi})}-1\right]
\mathscr{F}_{x}g(\xi)d\xi\right| \nonumber\\
&&\leq \left|\int_{|\xi|\leq \epsilon}e^{ix \xi}
\left[e^{it(\xi^{3}\pm \frac{1}{\xi})}-1\right]
\mathscr{F}_{x}g(\xi)d\xi\right|
\nonumber\\&&\qquad\qquad+\left|\int_{\epsilon
\leq |\xi|\leq 8}e^{ix \xi}\left[e^{it(\xi^{3}\pm
 \frac{1}{\xi})}-1\right]\mathscr{F}_{x}g(\xi)d\xi\right|\nonumber\\
&&\leq \int_{|\xi|\leq  \epsilon}|\mathscr{F}_{x}g(\xi)
|d\xi +C|t|\int_{ \epsilon\leq |\xi|\leq 8}
\frac{1}{|\xi|}|\mathscr{F}_{x}g(\xi)|d\xi\nonumber\\
&&\leq  \epsilon+\frac{C|t|}{\epsilon}\int_{|\xi|\leq8}|\mathscr{F}_{x}g(\xi)
|d\xi\leq C\left[\epsilon+\frac{|t|}{\epsilon}\right].\label{2.013}
\end{eqnarray}

This completes the proof of Lemma 2.4.

\begin{Lemma}\label{lem2.5}(Estimates  related to frequency-uniform decomposition)
Let   $|k|\leq8$ and  $g$ be
 a rapidly decreasing function. Then, $\forall \epsilon>0$, we have
\begin{align}
\left|U(t)\psi(D-k)g
-\psi(D-k)g\right|
\leq C\left[\epsilon+\frac{|t|}{\epsilon}\right].\label{2.014}
\end{align}
Here, $C$ is independent of $\epsilon, x,t$ and depends on $\|\mathscr{F}_{x}g\|_{L^{1}}$ .
\end{Lemma}
\noindent{\bf Proof.}   Since $g$ is a rapidly decreasing function, we have
\begin{eqnarray}
&&\int_{|\xi|\leq \epsilon}|\mathscr{F}_{x}g(\xi)|d\xi
\leq C\epsilon.\label{2.015}
\end{eqnarray}
For $\epsilon \leq |\xi|\leq 9$, we have
\begin{eqnarray}
&&\left|e^{it(\xi^{3}\pm \frac{1}{\xi})}-1\right|\leq
 |t|\left|\xi^{3}\pm \frac{1}{\xi}\right|\leq
 \frac{C|t|}{|\xi|}\leq \frac{C|t|}{\epsilon}\label{2.016}.
\end{eqnarray}
Thus, from (\ref{2.015})-(\ref{2.016}),  Since $g$ is a rapidly decreasing function, we have
\begin{eqnarray}
&&\left|U(t)\psi(D-k)g-\psi(D-k)g\right|\nonumber\\&&
=\left|\int_{\SR}e^{ix \xi}\left[e^{it(\xi^{3}\pm
\frac{1}{\xi})}-1\right]\psi(\xi-k)\mathscr{F}_{x}
g(\xi)d\xi\right| \nonumber\\
&&\leq \left|\int_{|\xi|\leq \epsilon}e^{ix \xi}
\left[e^{it(\xi^{3}\pm \frac{1}{\xi})}-1\right]
\psi(\xi-k)\mathscr{F}_{x}g(\xi)d\xi\right|
\nonumber\\&&\qquad\qquad+\left|\int_{|\xi|\geq
\epsilon}e^{ix \xi}\left[e^{it(\xi^{3}\pm \frac{1}
{\xi})}-1\right]\psi(\xi-k)\mathscr{F}_{x}g(\xi)d\xi\right|\nonumber\\
&&\leq2\int_{|\xi|\leq \epsilon}
|\mathscr{F}_{x}g(\xi)|d\xi +
2|t|\int_{\epsilon\leq |\xi|\leq 9}\frac{1}
{|\xi|}|\mathscr{F}_{x}g(\xi)|d\xi\nonumber\\
&&\leq C\epsilon+\frac{2|t|}{\epsilon}
\int_{\SR}|\mathscr{F}_{x}g(\xi)|d\xi\leq
 C\left[\epsilon+\frac{|t|}{\epsilon}\right].\label{2.017}
\end{eqnarray}

This completes the proof of Lemma 2.5.

\begin{Lemma}\label{lem2.6}
For $f\in L^{2}(\R)$, we have
\begin{eqnarray}
\left[\sum\limits_{k\in \z}|\psi(D-k)f|^{2}\right]^{\frac{1}{2}}\leq
\left\|f\right\|_{L^{2}(\SR)}\label{2.018}.
\end{eqnarray}
\end{Lemma}
\noindent{\bf Proof.} To obtain (\ref{2.018}), it suffices to prove
\begin{eqnarray}
\sum\limits_{k\in \z}|\psi(D-k)f|^{2}\leq
\left\|f\right\|_{L^{2}(\SR)}^{2}\label{2.019}.
\end{eqnarray}
By using the Cauchy-Schwarz inequality with respect to $\xi$, since $\supp\psi\subset [-1,1]$, we have
\begin{eqnarray}
&&\sum\limits_{k\in \z}|\psi(D-k)f|^{2}=\frac{1}{(2\pi)^{\frac{1}{2}}}\sum\limits_{k\in \z}\left|\int_{\SR}e^{i x\xi}\psi(\xi-k)\mathscr{F}_{x}f(\xi)d\xi\right|^{2}\nonumber\\
&&=\frac{1}{(2\pi)^{\frac{1}{2}}}\sum\limits_{k\in \z}\left|\int_{|\xi-k|\leq1}
e^{ix\xi}\psi(\xi-k)\mathscr{F}_{x}f(\xi)d\xi\right|^{2}\nonumber\\
&&\leq\sum\limits_{k\in \z}\left[\int_{|\xi-k|\leq1}|\psi(\xi-k)\mathscr{F}_{x}f(\xi)|^{2}d\xi
\int_{|\xi-k|\leq1}d\xi\right]\nonumber\\
&&\leq \sum\limits_{k\in \z}\int_{|\xi-k|\leq1}|\psi(\xi-k)\mathscr{F}_{x}f(\xi)|^{2}d\xi\nonumber\\
&&=\sum\limits_{k\in \z}\left\|\psi(\xi-k)\mathscr{F}_{x}f(\xi)\right\|_{L^{2}}^{2}.\label{2.020}
\end{eqnarray}
We claim
 \begin{eqnarray}
&& \|f\|_{L^{2}}^{2}\simeq\sum\limits_{k\in \z}\left\|\psi(\xi-k)\mathscr{F}_{x}f(\xi)\right\|_{L^{2}}^{2}\label{2.008}.
 \end{eqnarray}
 Now we prove the claim.
On one hand, from
\begin{eqnarray}
\mathscr{F}_{x}f(\xi)=\sum\limits_{k\in \z}\psi(\xi-k)\mathscr{F}_{x}f(\xi),\label{2.021}
\end{eqnarray}
by using the  Plancherel identity and $\supp\psi\subset [-1,1]$ and $\psi \geq0$,  we have
\begin{eqnarray}
&&\|f\|_{L^{2}}^{2}=\|\mathscr{F}_{x}f(\xi)\|_{L^{2}}^{2}=
\sum\limits_{k\in \z}\sum\limits_{l\in \z}\int_{\SR}\left[\psi(\xi-k)\mathscr{F}_{x}f(\xi)\right]\left[\psi(\xi-l)
\overline{\mathscr{F}_{x}f}(\xi)\right]d\xi\nonumber\\&&=
\sum\limits_{k\in \z}\int_{\SR}\left|\psi(\xi-k)\mathscr{F}_{x}f(\xi)\right|^{2}d\xi+\sum\limits_{k\in \z}\int_{\SR}
\psi(\xi-k)\mathscr{F}_{x}f(\xi)\psi(\xi-k+1)\overline{\mathscr{F}_{x}f}(\xi)d\xi\nonumber\\
&&\qquad+\sum\limits_{k\in \z}\int_{\SR}\psi(\xi-k)\mathscr{F}_{x}
f(\xi)\psi(\xi-k-1)\overline{\mathscr{F}_{x}f}(\xi)d\xi\nonumber\\
&&=\sum\limits_{k\in \z}\int_{\SR}\left|\psi(\xi-k)\mathscr{F}_{x}f(\xi)\right|^{2}d\xi+\sum\limits_{k\in \z}
\int_{\SR}\psi(\xi-k)\psi(\xi-k+1)|\mathscr{F}_{x}f(\xi)|^{2}d\xi\nonumber\\
&&\qquad+\sum\limits_{k\in \z}\int_{\SR}\psi(\xi-k)\psi(\xi-k-1)|\mathscr{F}_{x}
f(\xi)|^{2}d\xi\nonumber\\
&&\geq \sum\limits_{k\in \z}\int_{\SR}\left|\psi(\xi-k)\mathscr{F}_{x}f(\xi)\right|^{2}d\xi\label{2.022}.
\end{eqnarray}
 On the other hand, by using  the Cauchy-Schwarz inequality, we have
\begin{eqnarray}
&&\|f\|_{L^{2}}^{2}=\sum\limits_{k\in \z}\int_{\SR}\left|\psi(\xi-k)\mathscr{F}_{x}f(\xi)\right|^{2}d\xi+\sum\limits_{k\in \z}
\int_{\SR}\psi(\xi-k)\mathscr{F}_{x}f(\xi)\psi(\xi-k+1)\overline{\mathscr{F}_{x}f}(\xi)d\xi\nonumber\\
&&\qquad+\sum\limits_{k\in \z}\int_{\SR}\psi(\xi-k)\mathscr{F}_{x}f(\xi)\psi(\xi-k-1)\overline{\mathscr{F}_{x}f}(\xi)d\xi \nonumber\\
&&\leq \sum\limits_{k\in \z}\int_{\SR}\left|\psi(\xi-k)\mathscr{F}_{x}f(\xi)\right|^{2}d\xi+\sum\limits_{k\in \z}\left\|\psi(\xi-k)\mathscr{F}_{x}f(\xi)\right\|_{L^{2}}\left\|\psi(\xi-k+1)\overline{\mathscr{F}_{x}f(\xi)}\right\|_{L^{2}}\nonumber\\
&&\qquad+\sum\limits_{k\in \z}\left\|\psi(\xi-k)\mathscr{F}_{x}f(\xi)\right\|_{L^{2}}\left\|\psi(\xi-k-1)\overline{\mathscr{F}_{x}f(\xi)}\right\|_{L^{2}}\nonumber\\
&&\leq \sum\limits_{k\in \z}\left\|\psi(\xi-k)\mathscr{F}_{x}f(\xi)\right\|_{L^{2}}^{2}\nonumber\\&&\qquad+\left[\sum\limits_{k\in \z}\left\|\psi(\xi-k)\mathscr{F}_{x}f(\xi)\right\|_{L^{2}}^{2}\right]^{\frac{1}{2}}
\left[\sum\limits_{k\in \z}\left\|\psi(\xi-k+1)\overline{\mathscr{F}_{x}f(\xi)}\right\|_{L^{2}}^{2}\right]^{\frac{1}{2}}\nonumber\\
&&\qquad+\left[\sum\limits_{k\in \z}\left\|\psi(\xi-k)\mathscr{F}_{x}f(\xi)\right\|_{L^{2}}^{2}\right]^{\frac{1}{2}}
\left[\sum\limits_{k\in\z}\left\|\psi(\xi-k-1)\overline{\mathscr{F}_{x}f(\xi)}\right\|_{L^{2}}^{2}\right]^{\frac{1}{2}}\nonumber\\
&&=3\sum\limits_{k\in \z}\left\|\psi(\xi-k)\mathscr{F}_{x}f(\xi)\right\|_{L^{2}}^{2}\label{2.009}.
\end{eqnarray}
Combining (\ref{2.022}) with  (\ref{2.009}),  we have
\begin{eqnarray}
\sum\limits_{k\in \z}\left\|\psi(\xi-k)\mathscr{F}_{x}f(\xi)\right\|_{L^{2}}^{2}\leq \|f\|_{L^{2}}^{2}\leq 3\sum\limits_{k\in \z}\left\|\psi(\xi-k)\mathscr{F}_{x}f(\xi)\right\|_{L^{2}}^{2}. \label{2.0009}
\end{eqnarray}
Which implies the claim  (\ref{2.008}) holds.

Combining  (\ref{2.020}) with (\ref{2.008}),   we derive (\ref{2.019}).

This completes the proof of Lemma 2.6.

\begin{Lemma}\label{lem2.7}
For $f\in L^{2}(\R)$, we have
\begin{eqnarray}
\left[\sum\limits_{k\in \z}|\psi(D-k)U(t)f|^{2}\right]^{\frac{1}{2}}\leq
\left\|f\right\|_{L^{2}(\SR)}\label{2.023}.
\end{eqnarray}
\end{Lemma}
\noindent{\bf Proof.} To obtain (\ref{2.023}), it suffices to prove
\begin{eqnarray}
\sum\limits_{k\in \z}|\psi(D-k)U(t)f|^{2}\leq
\left\|f\right\|_{L^{2}}^{2}\label{2.024}.
\end{eqnarray}
By using the Cauchy-Schwarz inequality with respect to $\xi$,  since $\supp \psi \subset [-1,1]$, we have
\begin{eqnarray}
&&\sum\limits_{k\in \z}|\psi(D-k)U(t)f|^{2}=\frac{1}{(2\pi)^{\frac{1}{2}}}\sum\limits_{k\in \z}\left|\int_{\SR}e^{ix\xi}e^{it(\xi^{3}\pm\frac{1}{\xi})}\psi(\xi-k)\mathscr{F}_{x}f(\xi)d\xi\right|^{2}\nonumber\\
&&=\frac{1}{(2\pi)^{\frac{1}{2}}}\sum\limits_{k\in \z}\left|\int_{|\xi-k|\leq1}e^{ix\xi}e^{it(\xi^{3}\pm\frac{1}{\xi})}\psi(\xi-k)\mathscr{F}_{x}f(\xi)d\xi\right|^{2}\nonumber\\
&&\leq \sum\limits_{k\in \z}\left[\int_{|\xi-k|\leq1}|\psi(\xi-k)\mathscr{F}_{x}f(\xi)|^{2}d\xi\int_{|\xi-k|\leq1}d\xi\right]\nonumber\\
&&\leq \sum\limits_{k\in \z}\int_{|\xi-k|\leq1}|\psi(\xi-k)\mathscr{F}_{x}f(\xi)|^{2}d\xi\nonumber\\
&&=\sum\limits_{k\in \z}\left\|\psi(\xi-k)\mathscr{F}_{x}f(\xi)\right\|_{L^{2}}^{2}.\label{2.025}
\end{eqnarray}
Combining (\ref{2.022}) with (\ref{2.025}), we derive (\ref{2.024}).

This completes the proof of Lemma 2.7.

\bigskip
\bigskip

\noindent{\large\bf 3. Probabilistic estimates of some random series}
\setcounter{equation}{0}
\setcounter{Theorem}{0}

\setcounter{Lemma}{0}

\setcounter{section}{3}
In this section, we establish the probabilistic estimates of some random series. More precisely,
we use Lemmas 2.2, 2.4, 2.5,  3.1 to establish the probabilistic estimates of some random series which are just Lemmas 3.2-3.4 in this paper
 which play crucial role in establishing Theorem 1.3.

\begin{Lemma}\label{lem3.1}
Assume (\ref{1.06}). Then, there exists $C>0$ such that
\begin{eqnarray*}
\left\|\sum_{k\in\z}g_k(\omega)c_k\right\|_{L_{\omega}^p(\Omega)}\leq C\sqrt{p}\left\|c_k\right\|_{l^{2}(\z)}.
\end{eqnarray*}
for all $p\geq2 $
and $\{c_k\}\in l^{2}(\Z)$.
\end{Lemma}

For the proof of Lemma 3.1, we refer the readers to Lemma 3.1 of \cite{BTL}.

\begin{Lemma}\label{lem3.2}
 Let $g$ be  is a  rapidly  decreasing  function
and we denote by $g^{\omega}$ the randomization of $g$ as defined in (\ref{1.07}).
 Then,  for $\epsilon>0$ and $\alpha>0$, there exist  $C, C_{1}>0$ such that
\begin{eqnarray}
\mathbb{P}(\Omega_{1}^{c})\leq  C_{1} e^{-\left(\frac{\alpha}
{Ce\left[\epsilon+\frac{|t|}{\epsilon}\right]}\right)^{2}},\label{3.01}
\end{eqnarray}
where
$
\Omega_{1}^{c}=\left\{\omega\in \Omega:
\left|U(t)g^{\omega}-g^{\omega}\right|>\alpha\right\}.
$
Here, $C,C_{1}$ is independent of $\epsilon, x,t$.
\end{Lemma}
\noindent{\bf Proof.} Since $\left[P_{\geq8}+P_{\leq8}\right]g=g,$ we have
\begin{eqnarray}
&&\left\|U(t)g^{\omega}-g^{\omega}\right\|_{L_{\omega}^{p}(\Omega)}\leq I_{1}+I_{2},\label{3.02}
\end{eqnarray}
where
\begin{eqnarray}
I_{1}=\left\|U(t)P_{\geq8}g^{\omega}-P_{\geq8}g^{\omega}\right\|_{L_{\omega}^{p}(\Omega)},
I_{2}=\left\|U(t)P_{\leq8}g^{\omega}-P_{\leq8}g^{\omega}\right\|_{L_{\omega}^{p}(\Omega)}.\label{3.03}
\end{eqnarray}
Since $g$  is  a  rapidly  decreasing  function, by  using  Lemma 3.1,  we have
\begin{eqnarray}
&&I_{1}=\left\|U(t)P_{\geq8}g^{\omega}
-P_{\geq8}g^{\omega}\right\|_{L_{\omega}^{p}(\Omega)}\nonumber\\&&\leq C\sqrt{p}\left[\sum_{k\in \z}
\left|\int_{\SR}(e^{-it(\xi^{3}\pm\frac{1}{\xi})}-1)e^{ix \xi}\psi(\xi-k)\mathscr{F}
P_{\geq8}g(\xi)d\xi\right|^{2}\right]^{\frac{1}{2}}\nonumber\\
&&\leq C|t|\sqrt{p}\left[\sum_{k\in \z}\left|\int_{\SR}\left|\xi^{3}\pm \frac{1}{\xi}\right|\psi(\xi-k)
\mathscr{F}P_{\geq8}g(\xi)d\xi\right|^{2}\right]^{\frac{1}{2}}\nonumber\\
&&\leq C|t|\sqrt{p}\left[\sum_{k\in \z}\left|\int_{|\xi-k|\leq1}\left|\xi^{3}\pm \frac{1}{\xi}\right|\psi(\xi-k)
\mathscr{F}P_{\geq8}g(\xi)d\xi\right|^{2}\right]^{\frac{1}{2}}\nonumber\\
&&\leq C|t|\sqrt{p}\left[\sum_{k\in \z}\left[\int_{|\xi-k|\leq1}\left|\xi^{3}\pm \frac{1}{\xi}\right|^{2}|\psi(\xi-k)
\mathscr{F}P_{\geq8}g(\xi)|^{2}d\xi\int_{|\xi-k|\leq1}d\xi\right]\right]^{\frac{1}{2}}\nonumber\\
&&\leq C|t|\sqrt{p}\left[\sum_{k\in \z}\int_{\SR}|\xi|^{6}
|\psi(\xi-k)\mathscr{F}P_{\geq8}g(\xi)|^{2}d\xi\right]^{\frac{1}{2}}\nonumber\\
&&\leq C|t|\sqrt{p}\left[\sum_{k\in \z}\left\|\psi(D-k)P_{\geq8}g\right\|_{H^{3}}^{2}\right]^{\frac{1}{2}}\nonumber\\
&&\simeq |t|\sqrt{p}\left\|P_{\geq8}g\right\|_{H^{3}}\leq C|t|\sqrt{p}.\label{3.04}
\end{eqnarray}
Now we are going to justify the $\simeq$ appearing in the last line of (\ref{3.04}).
Here, by using the Cauchy-Schwarz inequality with respect to $\xi$ and (\ref{2.008}), we know
  \begin{eqnarray*}
  &&\left[\sum_{k\in \z}\left\|\psi(D-k)P_{\geq8}g(x)\right\|_{H^{3}}^{2}\right]^{\frac{1}{2}}
  =\left[\sum_{k\in \z}\left\|\psi(D-k)J^{3}P_{\geq8}g\right\|_{L^{2}}^{2}\right]^{\frac{1}{2}}\nonumber\\&&\simeq
  \left\|J^{3}P_{\geq8}g\right\|_{L^{2}}=\left\|P_{\geq8}g\right\|_{H^{3}},
  \end{eqnarray*}
where $J^{s}f=\mathscr{F}_{x}^{-1}\left(\langle\xi\rangle^{s}\mathscr{F}_{x}f(\xi)\right).$

From Lemmas 2.5 and 3.1, $\forall \epsilon>0$ such that $\epsilon<10^{-2}$, we have
\begin{eqnarray}
&&I_{2}= \left\|U(t)P_{\leq8}g^{\omega}
-P_{\leq8}g^{\omega}\right\|_{L_{\omega}^{p}(\Omega)}\nonumber\\&&\leq C\sqrt{p}
\left[\sum_{|k|\leq 10}\left|\int_{\SR}(e^{-it(\xi^{3}\pm\frac{1}{\xi})}-1)e^{ix \xi}\psi(\xi-k)\mathscr{F}P_{\leq8}g(\xi)d\xi\right|^{2}\right]^{\frac{1}{2}}\nonumber\\
&&\leq C\sqrt{p}\left[\sum_{|k|\leq 10}\left|\epsilon+\frac{|t|}{\epsilon}\right|^{2}\right]^{\frac{1}{2}}\leq C\sqrt{p}\left[\epsilon+\frac{|t|}{\epsilon}\right].\label{3.05}
\end{eqnarray}
From (\ref{3.02})-(\ref{3.05}),  $\forall \epsilon>0$ such that  $\epsilon<10^{-2}$, we have
\begin{eqnarray}
&&\left\|U(t)g^{\omega}-g^{\omega}\right\|_{L_{\omega}^{p}(\Omega)}\leq
C\sqrt{p}\left[\epsilon+\frac{|t|}{\epsilon}\right].\label{3.06}
\end{eqnarray}
Thus, from (\ref{3.06}), by using the Chebyshev inequality,  $\forall \epsilon>0$ such that $\epsilon<10^{-2}$ and $\forall\alpha>0$,   we have
\begin{eqnarray}
&&\mathbb{P}(\Omega_{1}^{c})\leq \int_{\Omega_{1}^{c}}\left[\frac{\left |U(t)g^{\omega}-g^{\omega}\right |}{\alpha}\right]^{p}d\mathbb{P}(\omega)\leq
\frac{\left\|U(t)g^{\omega}-g^{\omega}\right\|_{L_{\omega}^{p}}^{p}}{\alpha^{p}}\nonumber\\&&\leq
\left[\frac{C\sqrt{p}\left[\epsilon+\frac{|t|}{\epsilon}\right]}{\alpha}\right]^{p}.\label{3.07}
\end{eqnarray}
Take
\begin{eqnarray}
p=\left(\frac{\alpha}{Ce\left[\epsilon+\frac{|t|}{\epsilon}\right]}\right)^{2}\label{3.08}.
\end{eqnarray}
If $p\geq2$,
 we have
\begin{eqnarray}
\mathbb{P}(\Omega_{1}^{c})\leq e^{-p}=
e^{-\left(\frac{\alpha}{Ce\left[\epsilon+\frac{|t|}{\epsilon}\right]}\right)^{2}}.\label{3.09}
\end{eqnarray}
If $p\leq2,$
 we have
\begin{eqnarray}
&&\mathbb{P}\left(\Omega_{1}^{c}\right)\leq e^{2}e^{-2}\leq C_{1}e^{-\left(\frac{\alpha}{Ce\left[\epsilon+\frac{|t|}{\epsilon}\right]}\right)^{2}}\label{3.010}.
\end{eqnarray}
Here $C_{1}=e^{2}.$

This completes the proof of Lemma 3.2.

\begin{Lemma}\label{lem3.3}
 Let $h\in L^{2}(\R)$
 and we denote by $h^{\omega}$ the randomization of $h$ as defined in
 (\ref{1.07}).  Then, there exist  $C>0$ and $C_{1}>0$ such that
\begin{eqnarray}
\mathbb{P}\left(\Omega_{2}^{c}\right)\leq C_{1} e^{-\left(\frac{\alpha}
{Ce\|h\|_{L^{2}}}\right)^{2}}, \quad \text{for} ~\alpha>0,  \label{3.011}
\end{eqnarray}
where
\begin{eqnarray}
\Omega_{2}^{c}=\left\{\omega\in \Omega: \left|U(t)h^{\omega}\right|>\alpha\right\}.\label{3.012}
\end{eqnarray}
Here, $C,C_{1}$ is independent of $x,t$.
\end{Lemma}
\noindent {\bf Proof.} By using  Lemmas 3.1, 2.7,  we have
\begin{eqnarray}
&&\left\|U(t)h^{\omega}\right\|_{L_{\omega}^{p}(\Omega)}
=\left\|\sum\limits_{k\in \z}g_{k}(\omega)U(t)\psi(D-k)h\right\|_{L_{\omega}^{p}(\Omega)}
\nonumber\\
&&\leq C\sqrt{p}\left[\sum\limits_{k\in \z}|U(t)\psi(D-k)h|^{2}\right]^{\frac{1}{2}}\leq C\sqrt{p}\|h\|_{L^{2}}.\label{3.013}
\end{eqnarray}
Thus, by Chebyshev inequality, from (\ref{3.013}), we have
\begin{eqnarray}
&&\mathbb{P}\left(\Omega_{2}^{c}\right)\leq \int_{\Omega_{2}^{c}}\left[\frac{|U(t)h^{\omega}|}{\alpha}\right]^{p}d\mathbb{P}(\omega)
\leq \left(\frac{C\sqrt{p}\|h\|_{L^{2}}}{\alpha}\right)^{p}\label{3.023}.
\end{eqnarray}
Take
\begin{eqnarray}
p=\left(\frac{\alpha}{Ce\|h\|_{L^{2}}}\right)^{2}\label{3.024}.
\end{eqnarray}
If $p\geq2$,
 we have
\begin{eqnarray}
\mathbb{P}(\Omega_{2}^{c})\leq e^{-p}=
e^{-\left(\frac{\alpha}{Ce\|h\|_{L^{2}}}\right)^{2}}.\label{3.025}
\end{eqnarray}
If $p\leq2,$
  we have
\begin{eqnarray}
&&\mathbb{P}\left(\Omega_{2}^{c}\right)\leq e^{2}e^{-2}\leq C_{1}e^{-\left(\frac{\alpha}
{Ce\|h\|_{L^{2}}}\right)^{2}}\label{3.026}.
\end{eqnarray}
Here $C_{1}=e^{2}.$

This completes the proof of Lemma 3.3.

\begin{Lemma}\label{lem3.4}
Let $h\in L^{2}(\R)$
 and we denote by $h^{\omega}$ the randomization of $h$ as defined in
(\ref{1.07}). Then, there exist $C>0$ and $C_{1}>0$ such that
\begin{eqnarray}
\mathbb{P}(\Omega_{3}^{c})\leq C_{1}{\rm exp}\left[-\left(\frac{\alpha}
{Ce\|h\|_{L^{2}}}\right)\right]^{2},\quad \text{for} ~\alpha>0.\label{3.027}
\end{eqnarray}
Here
\begin{eqnarray*}
\Omega_{3}^{c}=\left\{\omega\in \Omega: |h^{\omega}|>\alpha\right\}
\end{eqnarray*}
and
 $C,C_{1}$ is independent of $x,t$.
\end{Lemma}
\noindent {\bf Proof.} By using the  Lemmas  3.1, 2.6,
  we have
\begin{eqnarray}
&&\hspace{-1cm}\left\|h^{\omega}\right\|_{L_{\omega}^{p}(\Omega)}
 =\left\|\sum\limits_{k\in \z}g_{k}(\omega)\psi(D-k)h\right\|_{L_{\omega}^{p}(\Omega)}
\nonumber\\
&&\leq C\sqrt{p}\left[\sum\limits_{k\in \z}|\psi(D-k)h|^{2}\right]^{\frac{1}{2}}\leq C\sqrt{p}\|h\|_{L^{2}}\label{3.028}.
\end{eqnarray}
Thus, by using the Chebyshev inequality,  from (\ref{3.028}),  we have
\begin{eqnarray}
&&\mathbb{P}\left(\Omega_{3}^{c}\right)\leq \int_{\Omega_{3}^{c}}\left[\frac{
|h^{\omega}|}{\alpha}\right]^{p}d\mathbb{P}(\omega)
\leq \left(\frac{C\sqrt{p}\|h\|_{L^{2}}}{\alpha}\right)^{p}\label{3.029}.
\end{eqnarray}
Take
\begin{eqnarray}
p=\left(\frac{\alpha}{Ce\|h\|_{L^{2}}}\right)^{2}\label{3.030}.
\end{eqnarray}
If $p\geq2$,
 we have
\begin{eqnarray}
\mathbb{P}(\Omega_{3}^{c})\leq e^{-p}=e^{-\left(\frac{\alpha
}{Ce\|h\|_{L^{2}}}\right)^{2}}.\label{3.031}
\end{eqnarray}
If $p\leq2,$
  we have
\begin{eqnarray}
&&\mathbb{P}\left(\Omega_{3}^{c}\right)\leq e^{2}e^{-2}\leq C_{1}e^{-\left(\frac{\alpha}
{Ce\|h\|_{L^{2}}}\right)^{2}}\label{3.032}.
\end{eqnarray}
Here $C_{1}=e^{2}.$

This completes the proof of Lemma 3.4.

\bigskip
\bigskip

\noindent{\large\bf 4. Proof of Theorem 1.1}
\setcounter{equation}{0}
\setcounter{Theorem}{0}

\setcounter{Lemma}{0}

\setcounter{section}{4}
In this section, we apply the density theorem and  Lemmas 2.1-2.4 to establish Theorem 1.1.

\noindent{\bf  Proof of Theorem 1.1.}
We firstly prove that if $f$ is rapidly decreasing function,
  \begin{eqnarray}
\left|U(t)f-f\right|\longrightarrow 0,\quad \text{as} ~t\longrightarrow0.\label{4.01}
\end{eqnarray}
From Lemmas 2.3, 2.4, $\forall \epsilon>0,$ we have
\begin{eqnarray}
\left|U(t)f-f\right|\leq C\left[\epsilon+\frac{|t|}{\epsilon}\right].\label{4.02}
\end{eqnarray}
When $|t|<\epsilon^{2}$, from (\ref{4.02}), we have
\begin{eqnarray}
\left|U(t)f-f\right|\leq 2C\epsilon.\label{4.03}
\end{eqnarray}
From (\ref{4.03}), we know that (\ref{4.01}) is valid.

\noindent When $f\in H^{s}(\R)(s\geq \frac{1}{4})$, by density theorem which can be seen
in Lemma 2.2 of \cite{Du}, for any $\epsilon>0$, there exists a rapidly decreasing function
 $g$ such that $f=g+h$, where $\|h\|_{H^{s}(\SR)}<\epsilon(s\geq \frac{1}{4}).$
Thus, we have
\begin{eqnarray}
\lim\limits_{t\longrightarrow0}\left|U(t)f-f\right|\leq \lim\limits_{t\longrightarrow0}\left|U(t)g-g\right|
+\lim\limits_{t\longrightarrow0}\left|U(t)h-h\right|
.\label{4.04}
\end{eqnarray}
For any $\alpha>0$ (fixed), we  define
\begin{eqnarray}
&&{\rm E_{\alpha}}=\left\{x\in \R: \lim\limits_{t\longrightarrow0}\left|U(t)f-f\right|>\alpha\right\}.\label{4.05}
\end{eqnarray}
Obviously, $E_{\alpha}\subset E_{1\alpha}\cup E_{2\alpha}$,
\begin{eqnarray}
&&{\rm E_{1\alpha}}=\left\{x\in \R: \lim\limits_{t\longrightarrow0}\left|U(t)
g-g\right|>\frac{\alpha}{2}\right\},\label{4.06}\\
&&{\rm E_{2\alpha}}=\left\{x\in \R: \lim\limits_{t\longrightarrow0}\left|U(t)h
-h\right|>\frac{\alpha}{2}\right\}.\label{4.07}
\end{eqnarray}
Obviously,
\begin{eqnarray}
E_{\alpha}\subset E_{1\alpha}\cup E_{2\alpha}.\label{4.08}
\end{eqnarray}
From Lemmas 2.3, 2.4,  we have
\begin{eqnarray}
\left|E_{1\alpha}\right|=0.\label{4.09}
\end{eqnarray}
Obviously,
\begin{eqnarray}
E_{2\alpha}\subset E_{21\alpha}\cup E_{22\alpha},\label{4.010}
\end{eqnarray}
where
\begin{eqnarray}
&&E_{21\alpha}=\left\{x\in \R: \sup\limits_{t>0}
\left|U(t)P_{\geq 8}h-P_{\geq8}h\right|>\frac{\alpha}{4}\right\},\label{4.011}\\
&&E_{22\alpha}=\left\{x\in \R:
\lim\limits_{t\longrightarrow0}\left|U(t)P_{\leq8}h-P_{\leq8}h\right|>\frac{\alpha}{4}\right\}\label{4.012}.
\end{eqnarray}
Thus, from  Lemma 2.1,  by  using  the Sobolev  embeddings  theorem $H^{\frac{1}{4}}(\R)\hookrightarrow L^{4}(\R)$,  we have
\begin{eqnarray}
&&\left| E_{21\alpha}\right|=\int_{E_{21\alpha}}dx\leq \int_{ E_{21\alpha}}
\frac{\left[\sup\limits_{t>0}\left|P_{\geq 8}U(t)h\right|\right]^{4}}
{\alpha_{1}^{4}}dx+\int_{ E_{21\alpha}}
\frac{\left|P_{\geq8}h\right|^{4}}{\alpha^{4}}dx\nonumber\\
&&\leq \frac{\left\|P_{\geq8}U(t)h\right\|_{L_{x}^{4}L_{t}^{\infty}}^{4}}{\alpha^{4}}
+\frac{\|P_{\geq8}h\|_{L_{x}^{4}}^{4}}{\alpha^{4}}\nonumber\\
&&\leq\frac{2C\left \|h\right\|_{H^{\frac{1}{4}}}^{4}}{\alpha^{4}}
\leq \frac{C\epsilon^{4}}{\alpha^{4}}\label{4.013}.
\end{eqnarray}
From Lemma 2.2 and $\epsilon$  is arbitrary, we have
\begin{eqnarray}
\left|E_{22\alpha}\right|=0.\label{4.014}
\end{eqnarray}
From  (\ref{4.09}), (\ref{4.013}) and (\ref{4.014}),  we have
\begin{eqnarray}
&&\left|E_{\alpha}\right|\leq\left|E_{1\alpha}\right|+\left|E_{2\alpha}\right|\leq
\left|E_{1\alpha}\right|+\left|E_{21\alpha}\right|+\left|E_{22\alpha}\right|\leq \frac{C\epsilon^{4}}{\alpha^{4}}\label{4.016}.
\end{eqnarray}
Thus, since  $\epsilon$  is arbitrary, from (\ref{4.016}),   we have
\begin{eqnarray}
\left|E_{\alpha}\right|=0.\label{4.017}
\end{eqnarray}
Thus, we have
\begin{eqnarray}
|U(t)f-f|\longrightarrow 0\label{4.018}
\end{eqnarray}
almost everywhere with respect to $x$ as $t$ goes to zero.

This completes the proof of Theorem 1.1.

\bigskip

\bigskip

\noindent {\large\bf 5. Proof of Theorem  1.2}

\setcounter{equation}{0}

 \setcounter{Theorem}{0}

\setcounter{Lemma}{0}

\setcounter{section}{5}

In this section, we  present the counterexample showing that $s\geq \frac{1}{4}$
is the necessary condition for the maximal function estimate related to free Ostrovsky equation.
 More precisely, we give the proof of Theorem 1.2.

 \noindent{\bf Proof of Theorem  1.2.}
 We  define $f_{k}=\frac{1}{\sqrt{2\pi}}\int_{\SR}e^{ix\xi}2^{-k(s+\frac{1}{2})}\chi_{2^{k}\leq |\xi|\leq 2^{k+1}}(\xi)d\xi$, obviously,
 \begin{eqnarray}
 \|f_{k}\|_{H^{s}}\simeq 1.\label{5.01}
 \end{eqnarray}
 Then, when $t\leq \frac{2^{-3k}}{100}$ and  $|x|\leq 2^{-k}$,
we have
\begin{eqnarray}
\|U(t)f_{k}\|_{L_{x}^{4}L_{t}^{\infty}}\gtrsim 2^{-k(s-1/4)}\label{5.02}.
\end{eqnarray}
Combining (\ref{5.01}) with (\ref{5.02}), for $s<\frac{1}{4}$, we have
\begin{align*}
\lim\limits_{k\rightarrow \infty}\frac{\left\|\sup\limits_{t>0}|U(t)f_{k}|\right\|_{L_{x}^{4}}}{\|f_{k}\|_{H^{s}(\SR)}}\geq C\lim\limits_{k\rightarrow \infty}\left\|\sup\limits_{t>0}|U(t)f_{k}|\right\|_{L_{x}^{4}}\geq C\lim\limits_{k\rightarrow \infty}2^{-k(s-1/4)}=+\infty.
\end{align*}
 From
\begin{eqnarray}
\|U(t)f_{k}\|_{L_{x}^{4}L_{t}^{\infty}}
\leq C\|f_{k}\|_{H^{s}}\label{5.03}
\end{eqnarray}
and (\ref{5.01})-(\ref{5.02}), we have
\begin{eqnarray}
2^{-k(s-\frac{1}{4})}\leq C\label{5.04}.
\end{eqnarray}
Hence, we know that
 for sufficiently large $k,$ when $s<\frac{1}{4}$, (\ref{5.04}) is invalid.

This completes the proof of Theorem 1.2.

\bigskip
\bigskip

\noindent {\large\bf 6. Proof of Theorem  1.3}

\setcounter{equation}{0}

 \setcounter{Theorem}{0}

\setcounter{Lemma}{0}

\setcounter{section}{6}
In this section, we apply Lemmas 3.2-3.4 and  the density theorem to prove Theorem 1.3.

\noindent{\bf Proof  of Theorem  1.3.}
We firstly prove that $\forall \alpha>0$(fixed) if $f$ is a rapidly decreasing function,
then
\begin{eqnarray}
\lim\limits_{t\longrightarrow 0}\mathbb{P}\left(\omega\in\Omega:|U(t)f^{\omega}(x)- f^{\omega}(x)|>\alpha\right)=0\label{6.01}
\end{eqnarray}
From Lemma 3.2,   $\forall  \epsilon>0$ such that 
$2Ce\epsilon(\ln\frac{3C_{1}}{\epsilon})^{\frac{1}{2}}<\alpha$,  we have
\begin{eqnarray}
\mathbb{P}\left(\left\{\omega\in \Omega:
\left|U(t)f^{\omega}-f^{\omega}\right|>\alpha\right\}\right)\leq
 C_{1}e^{-\left(\frac{\alpha}
{Ce\left[\epsilon+\frac{|t|}{\epsilon}\right]}\right)^{2}}\label{6.02}.
\end{eqnarray}
From (\ref{6.02}),  we know that when $|t|\leq  \epsilon^{2},$   we have
\begin{eqnarray}
\mathbb{P}\left(\left\{\omega\in \Omega:
\left|U(t)f^{\omega}-f^{\omega}\right|>\alpha\right\}\right)\leq
C_{1}e^{-\left(\frac{\alpha}{Ce\epsilon}\right)^{2}}\leq \epsilon\label{6.03}.
\end{eqnarray}
Hence, for any $\alpha>0$, we have
\begin{eqnarray}
\lim\limits_{t\longrightarrow 0}\mathbb{P}\left(\left\{\omega\in \Omega:
\left|U(t)f^{\omega}-f^{\omega}\right|>\alpha\right\}\right)=0.\label{1.025}
\end{eqnarray}
uniformly with respect to $x$.

Thus, we have proved (\ref{6.01}).

\noindent When $f\in L^{2}(\R)$, by density theorem which is Lemma 2.2
 in \cite{Du}, for any $\epsilon>0$, there exists a rapidly decreasing function
 $g$ such that $f=g+h$ which yields $f^{\omega}=g^{\omega}+h^{\omega}$, where $\|h\|_{L^{2}(\SR)}<\epsilon.$
Thus, we have
\begin{eqnarray}
\left|U(t)f^{\omega}-f^{\omega}\right|\leq \left|U(t)g^{\omega}-g^{\omega}\right|
+\left|U(t)h^{\omega}-h^{\omega}\right|
.\label{6.05}
\end{eqnarray}
From (\ref{6.05}), we have
\begin{eqnarray}
&&\left\{\omega\in \Omega: \left|U(t)f^{\omega}-f^{\omega}\right|>\alpha\right\}\nonumber\\&&\subset \left\{\omega\in \Omega: \left|U(t)
g^{\omega}-g^{\omega}\right|>\frac{\alpha}{2}\right\}\cup \left\{\omega\in \Omega: \left|U(t)h^{\omega}
-h^{\omega}\right|>\frac{\alpha}{2}\right\}.\label{6.06}
\end{eqnarray}
From Lemma 3.2, $\forall \epsilon>0$, when $|t|\leq \epsilon^{2}$, we have
\begin{eqnarray}
\mathbb{P}\left(\left\{\omega\in \Omega: \left|U(t)g^{\omega}-g^{\omega}\right|>\frac{\alpha}{2}\right\}\right)\leq
 C_{1}e^{-\left(\frac{\alpha}
{Ce\left[\epsilon+\frac{|t|}{\epsilon}\right]}\right)^{2}}\leq C_{1}e^{-\left(\frac{\alpha}{2Ce\epsilon}\right)^{2}}.\label{6.07}
\end{eqnarray}
From Lemmas 3.3, 3.4, we have
\begin{eqnarray}
&&\mathbb{P}\left(\left\{\omega\in \Omega: \left|U(t)h^{\omega}-h^{\omega}\right|>\frac{\alpha}{2}\right\}\right)\nonumber\\&&\leq
\mathbb{P}\left(\left\{\omega\in \Omega: \left|U(t)h^{\omega}\right|>\frac{\alpha}{4}\right\}\right)+\mathbb{P}\left(\left\{\omega\in \Omega: \left|h^{\omega}\right|>\frac{\alpha}{4}\right\}\right)\nonumber\\
&&\leq C_{1} e^{-\left(\frac{\alpha}
{Ce\|h\|_{L^{2}}}\right)^{2}}+C_{1} e^{-\left(\frac{\alpha}
{Ce\|h\|_{L^{2}}}\right)^{2}}=2C_{1} e^{-\left(\frac{\alpha}
{Ce\|h\|_{L^{2}}}\right)^{2}}.\label{6.08}
\end{eqnarray}
Combining  (\ref{6.06}), (\ref{6.07}) with (\ref{6.08}), we have that for any $\alpha>0$, $\forall\epsilon>0$ such that $2Ce\epsilon(\ln\frac{3C_{1}}{\epsilon})^{\frac{1}{2}}<\alpha$, taking $\|h\|_{L^{2}}<\epsilon$, the following inequality holds.
\begin{eqnarray}
&&\mathbb{P}\left(\left\{\omega\in \Omega: \left|U(t)f^{\omega}-f^{\omega}\right|>\alpha\right\}\right)\nonumber\\&&\leq \mathbb{P}\left(\left\{\omega\in \Omega: \left|U(t)g^{\omega}-g^{\omega}\right|>\alpha\right\}\right)+\left\{\omega\in \Omega: \left|U(t)h^{\omega}-h^{\omega}\right|>\alpha\right\}\nonumber\\&&\leq
C_{1}e^{-\left(\frac{\alpha}{2Ce\epsilon}\right)^{2}}+2C_{1} e^{-\left(\frac{\alpha}{Ce\|h\|_{L^{2}}}\right)^{2}}\nonumber\\&&\leq C_{1}e^{-\left(\frac{\alpha}{2Ce\epsilon}\right)^{2}}+2C_{1}e^{-\left(\frac{\alpha}{2Ce\epsilon}\right)^{2}}
=3C_{1}e^{-\left(\frac{\alpha}{2Ce\epsilon}\right)^{2}}\leq \epsilon.\label{6.09}
\end{eqnarray}
 From (\ref{6.09}),  we have
\begin{eqnarray}
\lim\limits_{t\longrightarrow 0}\mathbb{P}\left(\left\{\omega\in \Omega:
\left|U(t)f^{\omega}-f^{\omega}\right|>\alpha\right\}\right)=0\label{6.010}.
\end{eqnarray}

This completes the proof of Theorem 1.3.

\bigskip
\bigskip

\leftline{\large \bf Acknowledgments}
We are deeply indebted to the reviewers for his/her valuable suggestions which greatly improve the original version of our paper.
Wei Yan was supported by NSFC grants (No. 11771127) and the Young core Teachers program of Henan province under
grant number 2017GGJS044, Jiqiao Duan was supported by the NSF grant (No. 1620449) and NSFC grants (No. 11531006, 11771449) and Meihua Yang was supported by NSFC grants (No. 11971184).

\end{document}